\documentclass[final]{siamltex}
\usepackage{amssymb}
\usepackage{amsmath}

\newtheorem{remark}[theorem]{Remark}
\newtheorem{example}[theorem]{Example}

\def\norm#1{\|#1\|}

\def\span{{\mbox{\rm span}}}

\def\cF{\mathcal{F}}

\def\cH{\mathcal{H}}

\def\cW{\mathcal{W}}

\def\NN{\mathbb{N}}

\def\RR{\mathbb{R}}

\def\H{\mathcal{H}}
\def\<{\langle}
\def\>{\rangle}

\def\epsilon{\varepsilon}

\def\l{\lambda}
\def\t{\widetilde}

\newcommand{\ip}[2]{\langle#1,#2\rangle}

\title{Fusion Frames and Distributed Processing}

\author{Peter G. Casazza\thanks{Department of Mathematics,
University of Missouri, Columbia, MO 65211, USA ({\tt
pete@math.missouri.edu}). This author's research was
partially supported by NSF grant DMS-0405373.} \and
Gitta Kutyniok\thanks{Institute of Mathematics,
Justus--Liebig--University Giessen,
35392 Gies\-sen, Germany ({\tt gitta.kutyniok@math.uni-giessen.de}).
This author's research was partially supported by DFG
research fellowship KU 1446/5.}
\and Shidong Li\thanks{Department of Mathematics,
San Francisco State University,
San  Francisco,  CA   94132, USA ({\tt shidong@math.sfsu.edu}).
This author's research was  partially supported
by NSF RUI grant DMS-0406979.}}

\begin{document}

\maketitle

\begin{abstract}
Let $\{W_i\}_{i\in I}$ be a (redundant) sequence of subspaces each being endowed with a weight
$v_i$, and let $\cH$ be the closed linear span of the $W_i$'s, a composite Hilbert space.
Provided that $\{(W_i,v_i)\}_{i \in I}$ satisfies a certain property which controls the
weighted overlaps of the subspaces, it is called a {\em fusion frame}. These systems contain
conventional frames as a special case, however they go far ``beyond frame theory''.
In case each subspace $W_i$ is equipped with a frame system $\{f_{ij}\}_{j \in J_i}$ by which it
is spanned, we refer to $\{(W_i,v_i,\{f_{ij}\}_{j \in J_i})\}_{i \in I}$ as a {\em fusion
frame system}.

In this paper, we describe a weighted and distributed processing procedure that fuse together
information in all subspaces $W_i$ of a fusion frame system to obtain the global
information in $\cH$.  The weighted and distributed processing technique described
in fusion frames is not only a natural fit in distributed processing systems
such as sensor networks, but also an efficient scheme for parallel processing
of very large frame systems. We further provide an extensive study of the robustness
of fusion frame systems.
\end{abstract}



\pagestyle{myheadings}
\thispagestyle{plain}
\markboth{P.~G.~Casazza, G.~Kutyniok, and S.~Li}{Fusion Frames and Distributed Processing}


\section{Introduction}
Frames, which are systems that provide robust, stable and usually non-unique representations
of vectors, have been a focus of study in the last two decades in applications where
redundancy plays a vital and useful role, e.g., filter bank theory \cite{BHF98},
sigma-delta quantization \cite{BPY06}, signal and image processing \cite{CD02},
and wireless communications \cite{HP02}.

However, a number of new applications have emerged where the set-up can hardly be modeled
naturally by one single frame system. They generally share a common property that requires
distributed processing. Furthermore, we are often overwhelmed by a deluge of data assigned
to one single frame system, which becomes simply too large to be handled numerically. In these
cases it would be highly beneficial to split a large frame system into a set of (overlapping)
much smaller systems, and being able to process locally within each sub-system effectively.

A distributed frame theory relating to a set of local frame systems is clearly in demand.
In this paper we develop a suitable theory based on fusion frames, which  provides exactly
the framework not only to model these applications but also to provide efficient algorithms
with sufficient robustness.

\subsection{Applications under Distributed Processing Requirements}

A variety of applications require distributed processing.  Among them there are, for instance,
wireless sensor networks \cite{IB05}, geophones in geophysics measurements and studies \cite{CG06},
the physiological structure of ear and hearing systems \cite{RJ05b}.
To understand the nature, the constraints, and related problems of these applications, let us elaborate
a bit further on the example of wireless sensor networks.

In wireless sensor networks, sensors of limited capacity and power are spread in an area sometimes as
large as an entire forest to measure the temperature, sound, vibration, pressure, motion and/or pollutants.
In some applications, wireless sensors are placed in a geographical area to detect and characterize chemical,
biological, radiological, and nuclear material. Such a sensor system is typically redundant, and there is no
orthogonality among sensors, therefore each sensor functions as a frame element in the system.
Due to practical and cost reasons, most sensors employed in such applications have severe constraints in
their processing power and transmission bandwidth. They often have strictly metered power supply as well.
Consequently, a typical large sensor network necessarily divides the network into redundant
sub-networks -- forming a set of subspaces.
The primary goal is to have local measurements transmitted to a local sub-station within a subspace for
a subspace combining. An entire sensor system in such applications could have a number of such local
processing centers.  They function as relay stations, and have the gathered information further submitted
to a central processing station for final assembly.

In such applications, distributed/local processing is built in the problem formulation.
A staged processing structure is prescribed.  We will have to be able to process the
information stage by stage from local information and to eventually fuse them together
at the central station. We see therefore that a mechanism of coherently collecting
sub-station/subspace information is required.

Meantime, due to the often-time unpredictable nature of geographical factors, certain local
sensor systems are less reliable than others.   While facing the task of combining
local subspace information coherently, one has also reasons to consider weighting the more
reliable sets of substation information more than suspected less reliable ones. Consequently,
the coherent combination mechanism we just saw as necessary often requires a weighted structure
as well.  We will show that fusion frame systems are created to fit such weighted and coherent fusion needs.

\subsection{Parallel Processing of Large Frame Systems}
In case that a frame system is simply too large to handle effectively from
the numerical stand point, there are needs to divide the large system into smaller and parallel ones.
Like many parallel processing mechanisms, one may consider splitting the large
system into multiple small systems for simpler and parallel processing.   Evidently, the subdivision mechanism must
take into consideration a coherent combination after the subsystem processing.  To make the subdivision
mechanism more robust, one may not want to (sometimes it is also impossible to) split the large system in
an independent or orthogonal fashion.  Such a splitting and then a coherent combination must produce precisely
the original result if the system were to be processed globally.

Fusion frame systems are created to fit such needs as well.  Weighted coherent combination of subsystems (as
provided by fusion frame theory) is also useful in such applications where losses of some subsystem information
occur. Sometimes, weighted coherent combination is also useful from an efficient and approximation point of view.
Some approaches such as the domain decomposition method \cite{TW05} also use coherent combinations.  However
the fusion frame theory will provide a much more flexible framework that also takes local frames into account.

\subsection{Fusion frames}

In this article, we are interested in weighted sequences of subspaces with controlled
``overlaps''.  Each subspace is equipped with local frames aiming at the development of a framework for the
applications discussed above. In \cite{CK04}, two of the authors
studied redundant subspaces for the purpose of easing the construction of frames by building them
locally in (redundant) subspaces and then piecing the local frames together by employing
a special structure of the set of subspaces.  This was referred to as a {\em frame of subspaces}.
Related approaches were undertaken by Aldroubi, Cabrelli, and Molter \cite{ACM04a} and
Fornasier \cite{For03b}. A similar idea was also used by Aldroubi and Gr\"oching in a quite
different context in \cite{AG00}. Moreover, Bodmann, Kribs, and Paulsen  employed
Parseval frames of subspaces for optimal transmission of quantum states \cite{BKP06}.
Some further results on the theory from \cite{CK04} can be
found in \cite{AK05}, and an extension was derived by Sun \cite{Sun05,Sun06}, however
without any possibility of equipping the subspaces with an underlying structure.

We will employ some parts of this theory and set it into the context
studied in this paper. The structure of the overlapping subspaces, e.g., which relate to
the sub-networks of a wireless sensor network, will be modeled by employing the notion from \cite{CK04}.
To avoid confusion with the long existing term ``frame for subspaces'' and to emphasize the fact that this
mathematical object will provide a framework to fuse data in the subspaces, we decided to coin it {\em fusion
frame} in this context. In our situation it will become also essential to view a fusion
frame together with a set of local frames for its subspaces, in which case we will speak of a
{\em fusion frame system}.

We observe that fusion frames contain conventional frames as a special
case.  This theory goes thereby ``beyond frame theory''. It turns out that the fusion frame theory is
in fact much more delicate due to complicated relations between the structure of the sequence of
weighted subspaces and the local frames in the subspaces and due to the extreme sensitivity with
respect to changes of the weights.

Our main motivation is to study fusion frame systems
with respect to their reconstruction properties to not only provide a comprehensive model for
applications which require distributed processing and which employ a distributed structure due to
complexity reasons, but also to build efficient algorithms for fusion and reconstruction. We
provide a general reconstruction formula by employing a so-called fusion frame operator, derive
a variety of ways to fuse/reconstruct depending on the ability of the application to process
off-line or only in real time, and present an iterative algorithm. Since we are also concerned
with applications having the choice between distributed and centralized reconstruction we
further show that in very special cases those reconstructions are in fact performed by employing
the same set of vectors, thereby presenting situations where distributed reconstruction demonstrates
the same behavior as centralized reconstruction, e.g., with respect to noise.

As discussed above,  sensor networks in particular suffer significantly from disturbances of individual
sensors or even whole sub-networks in the form of, e.g., natural forces. This led us to the study of
stability of fusion frame systems not only under perturbations of the subspaces themselves, but even more
of the local frame vectors. In order to describe the properties of the affected sensor network
explicitly, we present several results which, in particular, give precise estimates for the changes
of certain properties of fusion frame systems.



\subsection{Contents}
The organization of this article is as follows.  In Section \ref{sec:fusion_frames}, the definition of
fusion frames and fusion frame systems and their fundamental characterization will be given. Examples
of fusion frames are presented, and connections of fusion frames with conventional frames will be discussed.
In Section \ref{sec:fusion_reconstruction}, several fusion frame reconstructions are presented.  These are
the coherent combinations we discussed earlier.  Both operator theoretical and its matrix representation
are considered.  An iterative fusion reconstruction is also constructed in this section.  Section
\ref{sec:robustness} is devoted to the robustness of fusion frames, in which the analysis of stability
of fusion frame systems to perturbations is extensively carried out.  Conclusion remarks and application
discussions are the subjects of the last section.


\section{Review of Frames and Notation}
\label{sec:frames}

A sequence $\cF=\{f_i\}_{i \in I}$ in a Hilbert space $\cH$ is a {\em frame} for $\cH$,
if there exist $0<A\leq B<\infty$ ({\em lower} and {\em upper frame bounds}) such that
\begin{equation} \label{eqn_frame}
A\|f\|^2 \leq \sum_{i \in I}|\< f, f_i\>|^2 \leq B\|f\|^2 \quad \mbox{for all } f \in \cH.
\end{equation}
The representation space associated with a frame is $\ell_2(I)$. In order to analyze
a signal $f \in \cH$, i.e., to map it into the representation space, the {\em analysis operator}
$T_\cF : \cH \to \ell_2(I)$ given by $T_\cF f = \{\ip{f}{f_i}\}_{i \in I}$ is applied.
The associated {\em synthesis operator}, which provides a mapping from the representation space
to $\cH$, is defined to be the adjoint operator $T_\cF^* : \ell_2(I) \to \cH$ which can
be computed to be $T_\cF^* (\{c_i\}_{i \in I}) = \sum_{i \in I} c_i f_i$. By composing
$T_\cF$ and $T_\cF^*$ we obtain the {\em frame operator}
\[ S_\cF : \cH \to \cH, \quad S_\cF f = T_\cF^* T_\cF f = \sum_{i \in I} \ip{f}{f_i} f_i.\]

Whenever $\cF=\{f_i\}_{i \in I}$ is a frame, we know that there exists at least one {\em dual
frame} $\{\tilde f_i\}_{i \in I}$ satisfying
\begin{equation} \label{eqn_frame_reconstruction}
     f=\sum_{i \in I}\< f, f_i\> \tilde f_i =\sum_{i \in I}\< f,\tilde f_i\> f_i \quad \mbox{for all } f \in \cH.
\end{equation}
When $\cF$ is a redundant (inexact) frame, there exist infinitely many dual frames $\{\tilde f_i\}_{i \in I}$
--~which can even be characterized \cite{Li95}~-- in which the {\em canonical dual frame} defined by
$\{S_\cF^{-1} f_i\}_{i \in I}$ is the one having the least square property
among all dual frames
$\{\tilde f_i\}_{i \in I}$.  That is, for all $f\in \cH$,
$\sum_{i \in I} |\< f, S^{-1}f_i\> |^2 \leq \sum_{i \in I} |\<f,\tilde f_i\>|^2$

Of particular interest are {\em $A$-tight frames}, i.e., if the frame bounds can be chosen as $A=B$ in
the frame definition (\ref{eqn_frame}). Provided (\ref{eqn_frame}) holds with $A=B=1$, we call $\cF$ a
{\em Parseval frame}. The advantage of working with these frames can be clearly seen by considering
the reconstruction formula (\ref{eqn_frame_reconstruction}). In these cases the canonical dual frame
equals $\{\frac{1}{A} f_i\}_{i \in I}$, and hence we obtain $f = \frac{1}{A} T_\cF^* T_\cF f$
for each $f \in \cH$, i.e., we can employ the frame elements for both the analysis and the
synthesis. There exist many procedures to construct tight or Parseval frames (cf. \cite{Cas04, CK07}).
However, Parseval frames with special properties are usually particularly difficult to construct,
see, e.g., \cite{SH03}.

For more details about the theory and applications of frames we refer the reader to
the books by Christensen \cite{Chr03}, Daubechies \cite{Dau92}, and Mallat \cite{Mal98}.


\section{Fusion Frames}
\label{sec:fusion_frames}

In this section, the notion of a {\em fusion frame} and a {\em fusion frame system} is introduced.
Discussions of the notion of redundancy for fusion frames is also provided.  We will put our focus
on the structure of the fusion frame operator and its connection with the fusion frame bounds, the
reason being that the fusion frame operator will become essential for studying distributed
fusion/reconstruction in Section \ref{sec:fusion_reconstruction}. Finally, the fact that our
theory goes ``beyond frame theory'' in the sense that conventional frames are a special case of fusion
frames will be discussed, and how much more sophisticated the theory of fusion frames turns
out to be will also be highlighted.


\subsection{Definition and Basic Properties}

We will start by stating the definition of a fusion frame.

\begin{definition}
Let $I$ be some index set, let $\{W_i\}_{i \in I}$ be a family of closed subspaces
in $\cH$, and let $\{v_i\}_{i \in I}$ be a family of weights, i.e., $v_i > 0$ for all $i \in I$.
Then $\{(W_i,v_i)\}_{i \in I}$ is a {\it fusion frame}, if there exist
constants $0 < C \le D < \infty$ such that
\begin{equation} \label{deffos}
C\|f\|^2 \le \sum_{i\in I} v_i^2 \|{\pi}_{W_i}(f)\|^2 \le D\|f\|^2
\ \ \mbox{for all $f\in \H$},
\end{equation}
where ${\pi}_{W_i}$ is the orthogonal projection onto the subspace $W_i$.
We call $C$ and $D$ the {\em fusion frame bounds}.
The family $\{(W_i,v_i)\}_{i\in I}$ is called a {\em $C$-tight fusion frame}, if in (\ref{deffos})
the constants $C$ and $D$ can be chosen so that $C=D$, a {\em Parseval fusion frame} provided that $C=D=1$
and an {\em orthonormal fusion basis} if $\H = \bigoplus_{i \in I} W_i$.
If $\{(W_i,v_i)\}_{i\in I}$ possesses an upper fusion frame bound, but not necessarily a lower bound,
we call it a {\it Bessel fusion sequence} with {\it Bessel fusion bound} $D$.
\end{definition}

Often it will become essential to consider a fusion frame together with a set
of local frames for its subspaces. In this case we will speak of a fusion frame system.

\begin{definition}
Let $\{(W_i,v_i)\}_{i \in I}$  be a fusion frame for $\cH$, and let $\{f_{ij}\}_{j \in J_i,\, i \in I}$
be a frame for $W_i$ for each $i \in I$. Then we call
$\{(W_i,v_i,\{f_{ij}\}_{j \in J_i})\}_{i \in I}$ a {\em fusion frame system} for $\cH$.
$C$ and $D$ are the associated {\em fusion frame bounds}, if they are the fusion frame bounds for
$\{(W_i,v_i)\}_{i \in I}$, and $A$ and $B$ are the {\em local frame bounds}, if these
are the common frame bounds for the {\em local frames} $\{f_{ij}\}_{j \in J_i}$ for each $i \in I$.
A collection of dual frames $\{\tilde f_{ij}\}_{j \in J_i}$, $i \in I$ associated with the local frames
will be called {\em local dual frames}.
\end{definition}

For a fusion frame system we have the following intriguing relation between properties
of the associated fusion frame and the sequence consisting of all local frame vectors,
in this sense it provides a link between local and global properties.
For the proof we refer to \cite[Thm. 3.2]{CK04}.

\begin{theorem} \label{theo:frame_local_global}
For each $i \in I$, let $v_i > 0$, let $W_i$ be a closed subspace
of $\cH$, and let $\{f_{ij}\}_{j \in J_i}$ be a frame for $W_i$
with frame bounds $A_i$ and $B_i$. Suppose
that
\[
0 < A = \inf_{i \in I} A_i \le \sup_{i \in I} B_i = B < \infty.
\]
Then the following conditions are equivalent.
\begin{enumerate}
\item $\{(W_i,v_i)\}_{i\in I}$ is a fusion frame for $\cH$.
\item $\{v_i f_{ij}\}_{j \in J_i,\, i \in I}$ is a frame for $\cH$.
\end{enumerate}
In particular, if $\{(W_i,v_i,\{f_{ij}\}_{j \in J_i})\}_{i\in I}$ is a fusion frame system for
$\cH$ with fusion frame bounds $C$ and $D$, then $\{v_i f_{ij}\}_{j \in J_i,\, i \in I}$ is a
frame for $\cH$ with frame bounds $AC$ and $BD$.
Also if $\{v_i f_{ij}\}_{i\in I,j \in J_i}$ is a frame for
$\cH$ with frame bounds $C$ and $D$, then $\{(W_i,v_i,\{f_{ij}\}_{j \in J_i})\}_{i\in I}$ is
a fusion frame system for $\cH$ with fusion frame bounds $\frac{C}{B}$ and $\frac{D}{A}$.
\end{theorem}

Tight frames play a vital role in frame theory due to the fact that they
provide easy reconstruction formulas, and also tight fusion frames will turn
out to be particularly useful for distributed reconstruction (cf. Section \ref{sec:fusion_reconstruction}).
The previous theorem
implies the following relation between tight fusion frames and tightness
of the collection of the local frames in a fusion frame system.

\begin{corollary} \label{coro:frame_local_global}
For each $i \in I$, let $v_i > 0$, let $W_i$ be a closed subspace
of $\cH$, and let $\{f_{ij}\}_{j \in J_i}$ be a Parseval frame for $W_i$.
Further, let $C$ be a constant. Then the following conditions are equivalent.
\begin{enumerate}
\item $\{(W_i,v_i)\}_{i\in I}$ is a $C$-tight fusion frame for $\cH$.
\item $\{v_i f_{ij}\}_{j \in J_i,\, i \in I}$ is a $C$-tight frame for $\cH$.
\end{enumerate}
\end{corollary}

By employing this result, the {\em redundancy} of a finite $C$-tight fusion
frame can be made precise in terms of the fusion frame bound.

\begin{proposition}
\label{prop:tight_redundancy} Let $\{(W_i,v_i)\}_{i=1}^n$ be a
$C$-tight fusion frame for $\cH$. Then we have
\[ C = \frac{\sum_{i=1}^n v_i^2 \dim W_i}{\dim \cH}.\]
\end{proposition}

\begin{proof}
Let $\{e_{ij}\}_{j =1}^{\dim W_i}$ be an orthonormal basis for $W_i$
for each $1 \le i \le n$. By Corollary \ref{coro:frame_local_global},
the sequence $\{v_i e_{ij}\}_{i=1,\,j =1}^{n, \quad \dim W_i}$ is a
$C$-tight frame for $\cH$. Employing \cite[Sec. 2.3]{CK03} yields that
\[ C = \frac{\sum_{i=1}^n\sum_{j=1}^{\dim W_i} \|v_i e_{ij}\|^2}{\dim \cH}
=\frac{\sum_{i=1}^n v_i^2 \dim W_i}{\dim \cH}.\]
\end{proof}

In this sense, we can interpret the frame bound $C$ as the redundancy of the tight
fusion frame $\{(W_i,v_i)\}_{i=1}^n$.

To enlighten the definitions let us consider the following example.

\begin{example}
\label{exa:finite_FF}
{\rm Since almost all applications require a finite model for their numerical treatment, we restrict
ourselfs to a finite dimensional space in this example.
Suppose $\{f_n\}_{n=1}^N$ is a frame for $\RR^M$ with frame bounds $A$, $B$.
Now we split $\{1,\ldots,N\}$ into $K$ sets $J_1,\ldots,J_K$, and define
$W_i = \span\{f_n\}_{n \in J_i}$, $1 \le i \le K$. Since in the finite-dimensional situation each
finite set of vectors forms a frame, in particular $\{f_n\}_{n \in J_i}$ is a frame for $W_i$
for each $1 \le i \le K$. Let $C$ and $D$ be a common lower and upper frame bound, respectively.
Theorem \ref{theo:frame_local_global} now implies that $\{(W_i,1,\{f_{n}\}_{n \in J_i})\}_{i=1}^K$
is a fusion frame system with fusion frame bounds $\frac{C}{B}$, $\frac{D}{A}$. Suppose that
by weighting the subspaces we can make $\{(W_i,v_i,\{f_{n}\}_{n \in J_i})\}_{i=1}^K$
 a tight frame, then the  fusion frame $\{(W_i,v_i)\}_{i=1}^K$ has redundancy
$(\sum_{i=1}^K v_i^2 \# J_i)/M$.}
\end{example}


\subsection{Fusion Frame Operator}

In frame theory an input signal is represented by a collection
of scalar coefficients that measure the projection of that signal onto each
frame vector. The representation space employed in this theory equals $\ell^2(I)$.
However, in fusion frame theory an input signal is represented by a collection
of {\em vector} coefficients that represent the projection (not just the projection
energy) onto each subspace. Therefore the representation space employed in this
setting is
\[\left ( \sum_{i\in I} \oplus W_i \right ) _{{\ell}^{2}} =
\{\{f_i \}_{i\in I}| f_i \in W_i  \mbox{ and } \{\norm{f_i}\}_{i\in
I} \in \ell^2(I)\}.\]

Let $\cW = \{(W_i,v_i) \}_{i\in I}$ be a fusion frame for $\H$. In order to map a signal
to the representation space, i.e., to analyze it, the {\em analysis operator} $T_\cW$
is employed, which is defined by
\[ T_\cW:\H \rightarrow \left(\sum_{i\in I}\oplus W_i \right)_{{\ell}_{2}}
\mbox{ with } \; T_\cW(f) = \{v_i{\pi}_{W_i}(f)\}_{i\in I}.\]
It can easily be shown that the {\em synthesis operator} $T_\cW^*$,
which is defined to be the adjoint operator, is given by
\[T_\cW^*: \left ( \sum_{i\in I} \oplus W_i \right )_{{\ell}_{2}} \rightarrow \H
\mbox{ with }\;  T_\cW(f) = \sum_{i\in I}v_i f_i, \;
f = \{f_i \}_{i\in I} \in \left(\sum_{i\in I}\oplus W_i \right)_{{\ell}_{2}}.\]
Now we can give the definition of a fusion frame operator.
The {\em fusion frame operator} $S_\cW$ for $\cW=\{(W_i,v_i)\}_{i\in I}$ is defined by
\[S_\cW(f) = T_\cW^*T_\cW(f)  = \sum_{i\in I}v_i^2{\pi}_{W_i}(f).\]

Interestingly, a fusion frame operator shows results similar to a frame operator
concerning invertibility. For the proof of the following result we refer to \cite[Prop. 3.16]{CK04}.

\begin{proposition} \label{prop:CI_S_DI}
Let $\{(W_i,v_i)\}_{i\in I}$ be a fusion frame for $\cH$ with fusion frame bounds $C$ and $D$.
Then the associated fusion frame operator $S_\cW$ is a positive and invertible
operator on $\H$ with $C\, $Id $\le S_\cW \le D\, $Id.
\end{proposition}

\subsubsection{Fusion frame operator in terms of local frames}

For the purpose of distributed fusion/reconstruction (see Section \ref{sec:fusion_reconstruction}),
employing a fusion frame system the fusion frame operator will indeed become essential. More precisely,
the inverse of the fusion frame operator will be employed. Therefore, a further investigation of the
fusion frame operator computationally is helpful.

We observe that the fusion frame operator can be expressed in terms of local frame operators as follows:

\begin{proposition}
Let $\{(W_i,v_i,\cF_i=\{f_{ij}\}_{j \in J_i})\}_{i\in I}$ be a fusion frame system for $\cH$, and let
$\t{\cF}_i=\{\tilde{f}_{ij}\}_{j\in J_i}$, $i \in I$ be associated local dual frames.
Then the associated fusion frame operator $S_\cW$ can be written as
\[
S_\cW =\sum_{i\in I} v_i^2 T^*_{\t{\cF}_i}T_{\cF_i}
        =\sum_{i\in I} v_i^2 T^*_{\cF_i}T_{\t{\cF}_i}.
\]
\end{proposition}

\begin{proof}
For all $f\in \cH$,
\[
 S_\cW f =\sum_{i\in I} v_i^2 \pi_{W_i} (f)
      =\sum_{i\in I}v_i^2\sum_{j\in J_i} \ip{f}{f_{ij}} \tilde f_{ij}
      = \sum_{i\in I}v_i^2\sum_{j\in J_i} \ip{f}{\tilde f_{ij}}  f_{ij}.
\]
Applying the definition of the analysis operators $T_{\cF_i}$, $T_{\t{\cF}_i}$
and the associated synthesis operators (see Section \ref{sec:frames}), the result follows immediately from here.
\end{proof}

\subsubsection{Matrix representation of the fusion frame operator}

For computational needs, let us further consider the fusion frame operator in {\em finite frame settings},
where the fusion frame operator will become the sum of (weighted) matrices of each subspace frame
operator (see also Example \ref{exa:finite_FF}).

Let $F_i$ be the frame matrices formed by frame vectors $\{f_{ij}\}_{j \in J_i}$ in the
column-by-column format
\[
   F_i \equiv \left( f_{i1} \, f_{i2}\, \cdots \, f_{i j_i}\right).
\]
Similarly, let $\tilde F_i$ be defined in the same way by the dual frame $\{\tilde f_{ij}\}_{j \in J_i}$.
Then the fusion frame operator associated with finite frames has the expression
\[
   S_\cW =\sum_{i\in I} v_i^2 F_i\tilde F_i^{H}
           =\sum_{i\in I} v_i^2 \tilde F_i F_i^{H},
\]
where $M^{H}$ stands for the {\em Hermitian transpose} of a matrix $M$.  Therefore, the evaluation of
the fusion frame operator and the inverse fusion frame operator in finite frame settings are quite
straightforward.  In practical applications, this will turn out to be very convenient.


\subsection{Fusion Frame Bounds}

Our first aim is to establish a connection between the fusion frame bounds and the
norm of the fusion frame operator. We achieve this by first showing that the boundedness of
the associated fusion frame operator is {\em equivalent} to the weighted sequence
of closed subspaces satisfying the fusion frame property.

\begin{proposition}
\label{prop:S_FF_boundedness}
Let $\{W_i\}_{i \in I}$ be closed subspaces in $\cH$, let $\{v_i\}_{i \in I}$ be
positive numbers, and let $S_\cW$ denote the fusion frame operator
associated with $\{(W_i,v_i)\}_{i \in I}$.
\begin{enumerate}
\item If $S_\cW \le D\, $Id, then $\{(W_i,v_i)\}_{i \in I}$ is a Bessel
fusion sequence with bound $D$.
\item If $S_\cW \ge C\, $Id, then $\{(W_i,v_i)\}_{i \in I}$ possesses the lower
fusion frame bound $C$.
\end{enumerate}
\end{proposition}

\begin{proof}
(i). Let $T_\cW$ denote the analysis operator associated with $\{(W_i,v_i)\}_{i \in I}$. Since
$S_\cW= T^*_\cW T_\cW$ and hence $\norm{T_\cW}^2 = \norm{S_\cW}$, for any $f \in \cH$ we obtain
\[ \sum_{i \in I} v_i^2 \norm{\pi_{W_i}(f)}^2
= \norm{T_\cW f}^2
\le \norm{T_\cW}^2 \norm{f}^2
= \norm{S_\cW} \norm{f}^2
\le D \norm{f}^2.\]

(ii). For all $f \in \cH$, we have
\[ \norm{T_\cW f}^2
=\ip{T^*_\cW T_\cW f}{f}
=\ip{S_\cW f}{f}
=\langle S^\frac12_\cW f,S^\frac12_\cW f\rangle
=\|S^\frac12_\cW f\|^2
\ge C\norm{f}^2.\]
\end{proof}

\begin{theorem}
\label{theo:S_FF}
Let $\{W_i\}_{i \in I}$ be closed subspaces in $\cH$, let $\{v_i\}_{i \in I}$ be
positive numbers, and let $S_\cW$ denote the fusion frame operator associated
with $\{(W_i,v_i)\}_{i \in I}$.
Then the following conditions are equivalent.
\begin{enumerate}
\item $\{(W_i,v_i)\}_{i \in I}$ is a fusion frame with fusion frame bounds $C$ and $D$.
\item We have $C\, $Id $\le S_\cW \le D\, $Id.
\end{enumerate}
Moreover, the fusion frame bounds are $\|S_\cW\|$ and $\|S^{-1}_\cW\|$.
\end{theorem}

\begin{proof}
(i) $\Rightarrow$ (ii). This is implied by Proposition \ref{prop:CI_S_DI}.

(ii) $\Rightarrow$ (i). This follows from Proposition \ref{prop:S_FF_boundedness}.
\end{proof}

Next we will study the behavior of the fusion frame operator under applying
a self-adjoint and invertible operator to the set of subspaces. In Subsection
\ref{subsec:beyond_frames}, this result will reveal essential differences
between fusion frame theory and classical frame theory.

\begin{proposition}
\label{prop:T_FF_SW}
Let $\{(W_i,v_i)\}_{i \in I}$ be a fusion frame for $\cH$ with associated fusion frame operator $S_\cW$,
and let $T$ be a self-adjoint and invertible operator on $\cH$. Then $\{(T W_i,v_i)\}_{i \in I}$ is
a fusion frame for $\cH$ with fusion frame operator $T S_\cW T^{-1}$.

In particular, $\{(S_\cW^{-1} W_i,v_i)\}_{i \in I}$ and $\{(S_\cW W_i,v_i)\}_{i \in I}$ both
possess the fusion frame operator $S_\cW$, and hence are fusion frames with the same fusion
frame bounds as $\cW$.
\end{proposition}

\begin{proof}
We recall the following basic fact from frame theory: Provided that $\cF=\{f_i\}_{i \in I}$ forms
a frame for a closed subspace $W$ of $\cH$ with frame operator $S_\cF$, then $\{Tf_{i}\}_{i \in I}$
is a frame for $TW$ with frame operator $TS_\cF T$, since
\[ \sum_{i \in I} \ip{f}{Tf_{i}}Tf_{i}
= T\left(\sum_{i \in I} \ip{Tf}{f_{i}}f_{i}\right)
= T S_\cF Tf.\]
Hence the dual frame of $\{Tf_{i}\}_{i \in I}$ is $\{T^{-1} S_\cF^{-1} f_{i}\}_{i \in I}$.

For each $i \in I$, let $\cF_i=\{f_{ij}\}_{j \in J_i}$ be a frame for $W_i$.
Then for all $f \in \cH$ we compute
\begin{eqnarray*}
\sum_{i \in I} v_i^2 \pi_{TW_i}(f)
& = & \sum_{i \in I} v_i^2  \left(\sum_{j \in J_i} \ip{f}{T^{-1}S_{\cF_i}^{-1}f_{ij}}Tf_{ij}\right)\\
& = & T \left(\sum_{i \in I} v_i^2 \pi_{W_i}(T^{-1}f)\right)\\
& = & T S_\cW T^{-1} f.
\end{eqnarray*}
The second part follows from here immediately by employing Theorem \ref{theo:S_FF}.
\end{proof}

This result shows in particular that the associated fusion frame is invariant under the
application of a self-adjoint and invertible operator, which commutes with the fusion frame
operator, to the set of subspaces.


\subsection{Beyond Frame Theory}
\label{subsec:beyond_frames}

Interestingly, frames can be shown to be a special case of fusion frames in a particular
sense with the natural meaning of the fusion frame bounds and the fusion frame operator.
We will make this precise in the following proposition.

\begin{proposition}
\label{prop:F_FF}
Let $\cF=\{f_i\}_{i \in I}$ be a frame for $\cH$ with frame bounds $A,B$.
Then $\{(\span\{f_i\},\norm{f_i})\}_{i \in I}$ is a fusion frame for $\cH$
with fusion frame bounds $A,B$ and fusion frame operator $S_\cF$.
\end{proposition}

\begin{proof}
Observe that for any $f \in \cH$,
\[\sum_{i \in I} \norm{f_i}^2 \pi_{\span\{f_i\}}(f)
=\sum_{i \in I} \norm{f_i}^2\ip{f}{\tfrac{f_i}{\norm{f_i}}} \tfrac{f_i}{\norm{f_i}}
=\sum_{i \in I}\ip{f}{f_i}f_i
=S_\cF f.\]
Hence the fusion frame operator for $\{(\span\{f_i\},\norm{f_i})\}_{i \in I}$ equals $S_\cF$.
Since $\{f_i\}_{i \in I}$ possesses the frame bounds $A$ and $B$, it follows
that $A\, $Id $\le S_\cF \le B\, $Id. Now Theorem \ref{theo:S_FF} implies
that $\{(\span\{f_i\},\norm{f_i})\}_{i \in I}$ is a fusion frame for $\cH$
with fusion frame bounds $A$ and $B$.
\end{proof}

This result seems to indicate that fusion frame theory is ``just'' a generalization of
frame theory. However, in the following remark we will enlighten the much more delicate
behavior of fusion frames. A variety of further essential differences will be revealed
by the results in the following sections.

\begin{remark}
{\rm To demonstrate the much more rich behavior of fusion frames in contrast to frames, we
consider a frame $\cF=\{f_i\}_{i \in I}$ for $\cH$ with frame bounds $A,B$.
Proposition \ref{prop:F_FF} implies that
\[\{(\span\{f_i\},\norm{f_i})\}_{i \in I}\]
is a fusion frame for $\cH$ with fusion frame bounds $A,B$ and fusion frame operator $S_\cF$.
Since $\{S_\cF^{-1}f_i\}_{i \in I}$ is the canonical dual frame for $\{f_i\}_{i \in I}$ with
bounds $B^{-1},A^{-1}$ and frame operator $S_\cF^{-1}$, also
\[\{(\span\{S_\cF^{-1}f_i\},\norm{S_\cF^{-1}f_i})\}_{i \in I}\]
is a fusion frame, but now with fusion frame bounds $B^{-1}, A^{-1}$ and fusion frame operator
$S_\cF^{-1}$.
Now it is possible to change the associated weights to ``move'' the associated fusion frame
operator back into the range of $A\, $Id and $B\, $Id. This is done by applying Proposition
\ref{prop:T_FF_SW} to $\{(\span\{f_i\},\norm{f_i})\}_{i \in I}$, which yields that
\[\{(S_\cF^{-1}\span\{f_i\},\norm{f_i})\}_{i \in I} = \{(\span\{S_\cF^{-1}f_i\},\norm{f_i})\}_{i \in I}\]
is a fusion frame with bounds $A,B$ and fusion frame operator $S_\cF$. This observation
not only reveals
 how much more sensitive fusion frames behave, but also indicates how critical
the selection of the weights can be.

Our observation is based on the fact that under the application of a self-adjoint and invertible operator
$T$ to both a frame and a fusion frame, the frame operator changes to $TS_\cF T$, however the fusion frame
operator changes to $T S_\cW T^{-1}$. Note that this fact for frames ensures that applying
$S_\cF^{-1}$ to the frame vectors yields a frame with frame operator $S_\cF^{-1}$. However, the analog
formula for fusion frames seems to make it almost impossible to construct a new fusion frame
having $S_\cW^{-1}$ as a fusion frame operator. We like therefore to state
this as an open question.}
\end{remark}


\section{Distributed Fusion/Reconstruction}
\label{sec:fusion_reconstruction}

Given a large set of data, some applications such as certain {\em data fusion problems} \cite{WaltzLlinas90}
require processing the data first locally by employing a frame structure, and then
fusing the (computed)
subspace information globally. This procedure is called {\em distributed fusion}, and obviously, the
second step can be modeled by employing the framework of fusion frame systems. If the initial data
comes from a decomposition of a signal with respect to a global frame such as in {\em sensor networks
problems} \cite{IB05}, and
the task consists in precisely reconstructing the initial signal via the procedure mentioned above,
we speak of {\em distributed reconstruction}. In this case we sometimes do have the choice of whether
either performing distributed or centralized reconstruction,
an issue that will be further elaborated in this section.

We shall first analyze the different distributed fusion procedures depending on whether it
is necessary to perform real time operations or whether it is possible to compute certain operations
off-line. We will then present an iterative algorithm for the computations in these procedures. Finally,
we discuss several aspects of distributed reconstruction versus centralized reconstruction, in particular
concerning the set of vectors employed to perform the reconstruction.


\subsection{Distributed Fusion Processing}

The first fundamental observation we make consists of the fact that
distributed fusion processing is
feasible in an elegant way by employing the inverse fusion frame operator.

\begin{proposition}
Let $\{(W_i,v_i)\}_{i\in I}$ be a fusion frame for $\cH$ with fusion frame operator $S_\cW$ and
fusion frame bounds $C$ and $D$. Then we have the reconstruction formula
\[ f = \sum_{i \in I} v_i^2 S_\cW^{-1} \pi_{W_i}(f) \quad \mbox{for all } f \in \H.\]
\end{proposition}

\begin{proof}
By Proposition \ref{prop:CI_S_DI}, for all $f \in \H$ we have
\[ f = S_\cW^{-1} S_\cW f =\sum_{i \in I} v_i^2 S_\cW^{-1} \pi_{W_i}(f).\]
\end{proof}

The fusion frame theory in fact provides two different approaches
for distributed fusion procedures.
For this, let $\{(W_i,v_i,\{f_{ij}\}_{j \in J_i})\}_{i\in I}$ be a fusion frame system for $\cH$,
and let $\{\tilde{f}_{ij}\}_{j\in J_i}$, $i \in I$ be associated local dual frames.

One distributed fusion procedure is from the local projections of each subspace:
\begin{equation}
f =\sum_{i\in I} v_i^2 S_\cW^{-1}\pi_{W_i}f
  =\sum_{i\in I}v_i^2 S_\cW^{-1}\left(\sum_{j\in J_i}\< f, f_{ij}\> \tilde f_{ij}\right)
  \quad \mbox{for all } f \in \H.
  \label{eqn_distReconst1}
\end{equation}
In this procedure, the local reconstruction takes place first in each subspace $W_i$, and
the inverse fusion frame is applied to each local reconstruction and combined together.

Another form of distributed fusion actually acts like a global reconstruction if the
coefficients of signal/function decompositions are available:
\begin{equation}
  f =\sum_{i\in I}v_i^2 \sum_{j\in J_i}\< f, f_{ij}\> \left(S_\cW^{-1}\tilde f_{ij}\right)
  \quad \mbox{for all } f \in \H.
  \label{eqn_distReconst2}
\end{equation}
The difference in this fusion procedure compared with global frame reconstruction lies in
the fact that the (global) dual frame $\{S_\cW^{-1}\tilde f_{ij}\}$ is first calculated at
the local level, and then fused into the global dual frame by applying the inverse fusion frame
operator.  This makes the evaluation of (global) duals much more efficient.

\begin{remark}
{\rm Depending on applications, some may require the fusion procedure via (\ref{eqn_distReconst1}) such
as in sensor networks \cite{IB05}, and geophones in geophysics measurements \cite{CG06},
whereas some may allow for fusion process via (\ref{eqn_distReconst2}) such as parallel processing
of large frame systems.  Let us examine the orders of
computation of the fusion procedures (\ref{eqn_distReconst1}) and (\ref{eqn_distReconst2}),
respectively.  Besides the operation of $S_\cW^{-1}$ in both equations, both fusion procedures
have the same number of multiplications.  However, (\ref{eqn_distReconst1}) typically has less
(but real time) inverse fusion frame operations.  Specifically, (\ref{eqn_distReconst1}) has $|I|$
operations of $S_\cW^{-1}$ over $|I|$ local reconstructions.   On the other hand, (\ref{eqn_distReconst2})
requires $\sum_{i\in I} |J_i|$ operations of $S_\cW^{-1}$ over local dual frames $\{f_{ij}\}_{j \in J_i,\, i \in I}$,
which is typically much larger than the $|I|$ operations in (\ref{eqn_distReconst1}).  It is nevertheless
equally important to point out that the much larger $S_\cW^{-1}$ operation requirement in
(\ref{eqn_distReconst2}) can be carried out ``off-line'', which often-times can be advantageous.}
\end{remark}


\subsection{Iterative Reconstruction}

Fusion frame reconstruction can be carried out iteratively as well, just like in frame
reconstructions \cite{Chr03}.   The specific mechanisms can also be divided in two
different ways,
depending on whether a local reconstruction actually takes place or not as given in
(\ref{eqn_distReconst1}) or (\ref{eqn_distReconst2}).

The first way we present refers to the distributed fusion procedure given by
(\ref{eqn_distReconst1}).

\begin{proposition} \label{prop_iterativeReconst}
Let $\{(W_i,v_i)\}_{i \in I}$ be a fusion frame in $\cH$  with fusion frame operator $S_\cW$
and fusion frame bounds $C$, $D$.
Further, let $f \in \cH$, and define the sequence $(f_n)_{n \in \NN_0}$ by
\[ f_n = \begin{cases}
0, & n=0, \\
f_{n-1} + \frac{2}{C+D}S_\cW(f-f_{n-1}), & n \ge 1.
\end{cases}\]
Then we have $f = \lim_{n \to \infty} f_n$ with the error estimate
\[ \norm{f-f_n} \le \left(\frac{D-C}{D+C}\right)^n \norm{f}.\]
\end{proposition}

\begin{proof}
Employing the construction of the sequence $(f_n)_{n \in \NN_0}$, for each $n \in \NN$
we obtain
\[ f-f_n = f-f_{n-1} - \frac{2}{C+D}S_\cW(f-f_{n-1}) = \left(I-\frac{2}{C+D}S_\cW\right)(f-f_{n-1}).\]
Interating this argument yields
\begin{equation} \label{eq:rec1}
f-f_n = \left(I-\frac{2}{C+D}S_\cW\right)^n(f-f_0).
\end{equation}
Now we use Proposition \ref{prop:CI_S_DI} to conclude that
\[ \left\langle\left(I-\frac{2}{C+D}S_\cW\right)f,f\right\rangle
\le \norm{f}^2 - \frac{2}{C+D}\ip{S_\cW f}{f}
\le \frac{D-C}{D+C}\norm{f}^2.\]
In a similar way, we can show that
\[ \left\langle\left(I-\frac{2}{C+D}S_\cW\right)f,f\right\rangle
\ge - \frac{D-C}{D+C}\norm{f}^2.\]
Combining these two estimates yields
\[\norm{I-\frac{2}{C+D}S_\cW} \le \frac{D-C}{D+C} < 1.\]
Applying this estimate to \eqref{eq:rec1}, we finally obtain
\[ \norm{f-f_n}_2 \le \norm{I-\frac{2}{C+D}S_\cW}^n \norm{f-f_0}
\le \left(\frac{D-C}{D+C}\right)^n \norm{f}.\]
\end{proof}

Thus every $f \in \cH$ can be reconstructed from the fusion frame coefficients $T_\cW(f) =
\{v_i \pi_{W_i}(f)\}_{i \in I}$, since $S_\cW f$ does only require the knowledge of those
and of the sequence of weights $\{v_i\}_{i \in I}$. In each iteration, $S_\cW f$ is always known.
The only significant computation is $S_\cW f_n$, which can be easily carried out using the local
frame structure of each subspace.

We remark that an application of the Chebyshev method or the conjugate gradient
method as done by Gr\"ochenig \cite{Gro93} for the frame algorithm should lead to faster
convergence.

\medskip

Finally, we discuss an interactive way to compute the distributed fusion procedure given by
(\ref{eqn_distReconst2}).

\begin{remark}
{\rm In some applications,
if the local measurements/local frame coefficients are preserved,
the final reconstruction can also be done in a ``global'' fashion with distributed evaluation
of (global) duals through local dual frames $\{\tilde f_{ij}\}_{j \in J_i}$, $i \in I$
associated with a fusion frame system $\{(W_i,v_i,\{f_{ij}\}_{j \in J_i})\}_{i\in I}$.

This follows from equation (\ref{eqn_distReconst2}), in which the evaluation of local duals
$\{\tilde f_{ij}\}_{j \in J_i}$ is carried out distributively, and the evaluation of
$\{S_\cW^{-1}\tilde f_{ij}\}_{j \in J_i,\, i \in I}$ can be carried out iteratively as described in
Proposition \ref{prop_iterativeReconst}.}
\end{remark}


\subsection{Distributed Reconstruction and (Global) Dual Frames}

The purpose of this section is to study the sequence of vectors employed for distributed
reconstruction and to compare distributed with centralized reconstruction. For this, let
$\{(W_i,v_i,\{f_{ij}\}_{j \in J_i})\}_{i\in I}$ be a fusion frame system for $\cH$
with local frame bounds $A$, $B$,
and let $\{\tilde{f}_{ij}\}_{j\in J_i}$, $i \in I$ be associated local dual frames.
Since, by Theorem~\ref{theo:frame_local_global}, the sequence $\cF=\{v_i f_{ij}\}_{j \in J_i,\, i \in I}$ is a
frame for $\cH$, we might
consider the situation that we are given the (global) frame coefficients $\{\< f, v_if_{ij}\>\}_{j \in J_i,\, i \in I}$
of a signal $f \in \cH$. For some applications, which do not enforce distributed reconstruction, we might have
two ways to reconstruct $f$. The (global) dual frame $\{S_\cF^{-1} v_i  f_{ij}\}_{j \in J_i,\, i \in I}$ could
be used to perform {\em centralized reconstruction}, i.e., to compute
\[f =\sum_{i\in I} \sum_{j\in J_i}\< f, v_if_{ij}\> \left(S_\cF^{-1} v_i  f_{ij}\right).\]
Or, in order to reduce the complexity, we might employ the associated fusion frame operator $S_\cW$ to
perform {\em distributed reconstruction}, and obtain (compare \eqref{eqn_distReconst2})
\[f =\sum_{i\in I} \sum_{j\in J_i}\< f, v_if_{ij}\> \left(S_\cW^{-1}v_i\tilde f_{ij}\right).\]
In the sequel we will discuss the difference between the sequences
$\{S_\cF^{-1} v_i  f_{ij}\}_{j \in J_i,\, i \in I}$
and $\{S_\cW^{-1} v_i \tilde{f}_{ij}\}_{j \in J_i,\, i \in I}$ in more detail.

Our first result shows that indeed $\{S_\cW^{-1} v_i \tilde{f}_{ij}\}_{j \in J_i,\, i \in I}$ is
{\em a dual frame} for $\cF$, but not necessarily the {\em canonical dual frame}.

\begin{proposition} \label{prop_dualFF}
Let $\{(W_i,v_i,\{f_{ij}\}_{j \in J_i})\}_{i\in I}$ be a fusion frame
system for $\cH$ with associated fusion frame operator $S_\cW$, local frame bounds $A,B$ and
local dual frames $\{\tilde{f}_{ij}\}_{j \in J_i}$, $i \in I$.
Then $\{S_\cW^{-1} v_i \tilde{f}_{ij}\}_{j \in J_i,\, i \in I}$
is a dual frame for the frame $\{v_i f_{ij}\}_{j \in J_i,\, i \in I}$.
\end{proposition}

\begin{proof}
First we note that $\{v_i f_{ij}\}_{j \in J_i,\, i \in I}$ indeed forms a frame by
Theorem \ref{theo:frame_local_global}. Employing the self-adjointness of $S_\cW$
(Proposition \ref{prop:CI_S_DI}), we have for all $f \in \cH$,
\begin{eqnarray*}
\sum_{i \in I} \sum_{j \in J_i} \ip{f}{S_\cW^{-1} v_i \tilde{f}_{ij}} v_i f_{ij}
& = & \sum_{i \in I} v_i^2 \sum_{j \in J_i} \ip{S_\cW^{-1} f}{\tilde{f}_{ij}} f_{ij}\\
& = & \sum_{i \in I} v_i^2 \sum_{j \in J_i} \ip{\pi_{W_i} (S_\cW^{-1} f)}{\tilde{f}_{ij}} f_{ij}\\
& = & \sum_{i \in I} v_i^2 \pi_{W_i}(S_\cW^{-1} f)\\
& = &  S_\cW (S_\cW^{-1} f)\\
& = & f,
\end{eqnarray*}
which finishes the proof.
\end{proof}

It is interesting to observe that a ``dual'' relation also holds. We wish to mention that this
property does not have quite the same correspondence in conventional frames as well.

\begin{proposition}
Let $\{(W_i,v_i,\{f_{ij}\}_{j \in J_i})\}_{i\in I}$ be a fusion frame
system for $\cH$ with associated fusion frame operator $S_\cW$, local frame bounds $A,B$ and
local dual frames $\{\tilde{f}_{ij}\}_{j \in J_i}$, $i \in I$.
Then $\{v_i \tilde f_{ij}\}_{i\in I, j\in J_i}$ is a frame for $\cH$
and $\{S_\cW^{-1} v_i f_{ij}\}_{i\in I, j\in J_i}$ is a dual frame for it.
\end{proposition}

\begin{proof}
The fact that $\{v_i \tilde f_{ij}\}_{i\in I, j\in J_i}$ is a frame for $\cH$
follows again from Theorem \ref{theo:frame_local_global}, since each local
frame $\{f_{ij}\}_{j \in J_i}$, $i \in I$ has frame bounds $A^{-1}$, $B^{-1}$
and thus possess a common lower and upper bound.

Now using the fact that $\{v_i f_{ij}\}_{i\in I, j\in J_i}$ and
$\{S_\cW^{-1} v_i \tilde f_{ij}\}_{i\in I, j\in J_i}$ are a pair of dual frames of $\cH$
by Proposition \ref{prop_dualFF}, we have for all $f\in\cH$,
{\allowdisplaybreaks \begin{eqnarray*}
 f &=& \sum_{i\in I}\sum_{j\in J_i}\< f, v_i f_{ij}\> S_\cW^{-1}v_i
      \tilde f_{ij} \\
  &=& S_\cW^{-1}\left(\sum_{i\in I}\sum_{j\in J_i}\< f, v_i f_{ij}\>
      v_i \tilde f_{ij}\right) \\
  &=& S_\cW^{-1}\left(\sum_{i\in I} v_i^2\sum_{j\in J_i}
      \< f, \tilde f_{ij}\> f_{ij}\right)\\
  &=& \sum_{i\in I}\sum_{j\in J_i}\< f, v_i\tilde f_{ij}\>
      S_\cW^{-1}v_i f_{ij}.
\end{eqnarray*}}
\end{proof}

In order to compare distributed reconstruction with centralized reconstruction, it is
essential to understand when $\{S_\cW^{-1} v_i \tilde{f}_{ij}\}_{j \in J_i,\, i \in I}$ equals
the canonical dual frame of the frame $\cF=\{v_i f_{ij}\}_{j \in J_i,\, i \in I}$ (compare also
\eqref{eqn_distReconst1}),
since in these particular cases distributed and centralized reconstruction coincide. In general,
this certainly need not be the case due to the observation that if we have a Parseval fusion frame,
then
\[ S_\cF = \sum_{i \in I} S_{\cF_i} \pi_{W_i},\]
with the $S_{\cF_i}$ being the local frame operators, and hence due to the occurring cross terms,
\[ \{v_i S_{\cF_i}^{-1} f_{ij}\}_{j \in J_i,\, i \in I}
\neq  \{v_i S_\cF^{-1} f_{ij}\}_{j \in J_i,\, i \in I}.\]
However, the following results give some special cases in which distributed and centralized reconstruction
indeed coincide.

\begin{proposition}
\label{prop:centr=distr_duals}
Let $\{(W_i,v_i,\cF_i=\{f_{ij}\}_{j \in J_i})\}_{i\in I}$ be a fusion frame
system for $\cH$ with associated fusion frame operator $S_\cW$, local frame bounds $A,B$ and
local dual frames $\{\tilde{f}_{ij}\}_{j \in J_i}$, $i \in I$.
If $\{(W_i,v_i)\}_{i \in I}$ is an orthogonal fusion basis or $\{f_{ij}\}_{j \in J_i}$
is a Parseval frame sequence for all $i \in I$, then $\{S_\cW^{-1} v_i S_{\cF_i}^{-1} f_{ij}\}_{j \in J_i,\, i \in I}$
is the canonical dual frame of the frame $\{v_i f_{ij}\}_{j \in J_i,\, i \in I}$.
\end{proposition}

\begin{proof}
By Theorem \ref{theo:frame_local_global}, the sequence $\cF=\{v_i f_{ij}\}_{j \in J_i,\, i \in I}$
forms a frame for $\cH$, and we denote its frame operator by $S_\cF$.
If $\{(W_i,v_i)\}_{i \in I}$ is an orthogonal fusion basis, then $S_\cW = $ Id and $S_\cF =
\sum_{i \in I} \oplus S_{\cF_i} \pi_{W_i}$, and hence $S_\cF^{-1} = \sum_{i \in I} \oplus S_{\cF_i}^{-1} \pi_{W_i}$.
Provided that $\{f_{ij}\}_{j \in J_i}$ is a Parseval frame sequence for all $i \in I$, we
have $S_{\cF_i} = $ Id for all $i \in I$, and we further obtain
for all $f \in \cH$,
\[S_\cW f = \sum_{i \in I} v_i^2 \pi_{W_i}(f) =  \sum_{i \in I} \sum_{j \in J_i}
\ip{f}{v_i f_{ij}}v_i f_{ij} = S_{\cF} f.\]
In both cases the claim follows immediately from here.
\end{proof}

Finally, we would like to point out a surprising fact concerning the situation of having
Parseval frames spanning the subspaces of a fusion frame, which arises from this result.

\begin{remark}
{\rm Let $\{(W_i,v_i,\cF_i=\{f_{ij}\}_{j \in J_i})\}_{i\in I}$ be a fusion frame
system for $\cH$ with associated fusion frame operator $S_\cW$, and let $\cF_i$ be
Parseval frames for all $i \in I$. By the previous result, the operator $S_\cW$ is
{\em independent}
of the choice of the Parseval frame, since $S_\cW$ always equals the frame operator
of the frame $\{v_i f_{ij}\}_{j \in J_i,\, i \in I}$. The intuitive reason for this is that
provided we take Parseval frames for the subspaces, the frame property of the total
collection of frame elements completely mirrors the behavior of the fusion frame.}
\end{remark}

Furthermore, we would like to briefly mention the impact of Proposition \ref{prop:centr=distr_duals}
on the noise reduction under distributed and centralized reconstruction of this result.

\begin{remark}
{\rm In \cite{RJ05a} Rozell and Johnson studied noise reduction under
distributed reconstruction versus centralized reconstruction in finite dimensional
Hilbert spaces. They used additive zero mean white noise, which they added to all
frame coefficients with respect to the frame $\cF = \{v_i f_{ij}\}_{j \in J_i,\, i \in I}$
and derived bounds for the mean square error for distributed and centralized reconstruction.
By numerical simulation they showed that randomly adding elements to a fusion frame not
only improves the mean square error for distributed reconstruction, but even forces it
to converge to the mean square error for centralized reconstruction. Goyal,
Vetterli, and Thao \cite{GVT98} proved that a frame becomes asymptotically tight
if elements are randomly added. By employing this result it was argued in \cite{RJ05a}
that some very restrictive conditions on the local frames and the fusion fame structure
might be fulfilled in the limit.

Applying Proposition \ref{prop:centr=distr_duals} now gives a broader picture for this
intriguing phenomenon. By \cite{GVT98},  the local frames become asymptotically
tight with the same frame bound, say $A$, as more and more random elements are added.
Hence, by Proposition \ref{prop:centr=distr_duals},
the sequence $\{S_\cW^{-1} v_i S_{\cF_i}^{-1} \frac{1}{\sqrt{A}}f_{ij}\}_{j \in J_i,\, i \in I}$,
which is employed for the distributed reconstruction as outlined in \eqref{eqn_distReconst1}
(except the constant $\frac{1}{\sqrt{A}}$), converges to the canonical dual frame of the frame
$\{v_i f_{ij}\}_{j \in J_i,\, i \in I}$, which is used for the centralized reconstruction.
Taking the fact into account that the additional multiplicative constant does not play any role
concerning the noise reduction ability, reveals a reason for the phenomenon described in
\cite{RJ05a} and does not require restrictive conditions.}
\end{remark}


\section{Robustness of Fusion Frame Systems}
\label{sec:robustness}

In this section we analyze the stability of fusion frame systems under perturbations
of both the subspaces which constitute a fusion frame and the local frame
vectors contained in the subspaces. The reason for this is that on the one
hand, for instance, several complete groups of geophones \cite{CG06} might be moved
to a slightly different location to adjust for transmission conditions,
and on the other hand, for instance, in wireless sensor networks the location
of single sensors might be changed slightly due to the impact of natural forces
\cite{IB05}.
Therefore it is essential to study the robustness of fusion frame systems under these
two different impacts.

Thus, with these practical aspects in mind, we proceed by first examining
perturbations of the subspaces in Section \ref{subsec:pert_subspaces} and
secondly studying robustness of a fusion frame system under perturbations of the
associated local frames in Section \ref{subsec:pert_frames}.

\subsection{Perturbation of the Fusion Frame}
\label{subsec:pert_subspaces}

First we would like to point out one fundamental problem with
perturbations of fusion frames which is the cause of the serious technicalities
in these results. Since the main ingredient in the definition of a fusion frame are
the orthogonal projections onto a set of subspaces, it would be natural to consider
perturbations of these projections. However, there is no such thing as a
perturbation of a projection.  This means, that if $P$ and $Q$ are projections
on $\cH$, $0\le \l_1,\l_2 < 1$, and
$$
\|Pf - Qf\| \le \l_1 \|Pf\| + \l_2 \|Qf\| \quad \mbox{for all $f\in \cH$},
$$
then it follows that $P=Q$.  This can be easily seen by way of contradiction
as follows: If $P\not= Q$, then there exists a vector $f\in \cH$ so that $f\perp P(\cH)$,
but also satisfying $Qf \not= 0$ (or vice-versa).  This yields
$$
\|Pf-Qf\| = \|Qf\| \le \l_1 \|Pf\| + \l_2 \|Qf\| = \l_2 \|Qf\|,
$$
which is a contradiction.


Therefore, we define $(\l_1,\l_2)$-perturbations of sequences by employing
the canonical Paley-Wiener-type definition:

\begin{definition}
Let $\{W_i\}_{i \in I}$ and $\{\t{W}_i\}_{i \in I}$ be closed subspaces in
$\cH$, let $\{v_i\}_{i \in I}$ be positive numbers, and let
$0 \le \lambda_1,\lambda_2 <1$ and $\epsilon > 0$. If
\[
\|(\pi_{W_i}-\pi_{\t{W}_i})f\| \le \lambda_1 \|\pi_{W_i}f\| + \lambda_2
\|\pi_{\t{W}_i}f\|
+ \epsilon \|f\| \quad \mbox{for all } f\in \cH \mbox{ and } i \in I,
\]
then we say that $\{(\t{W}_i,v_i)\}_{i \in I}$ is a
{\em $(\lambda_1,\lambda_2,\epsilon)$-perturbation} of $\{(W_i,v_i)\}_{i \in I}$.
\end{definition}

Employing this definition, we derive the following result about robustness of
fusion frames under small perturbations of the associated subspaces.

\begin{proposition}\label{propp}
Let $\{(W_i,v_i)\}_{i \in I}$ be a fusion frame for $\cH$ with bounds $C$, $D$. Choose
$0 \le \l_1 < 1$ and $\epsilon >0$ such that
\[ (1- \l_1)\sqrt{C} - \epsilon \left ( \sum_{i\in I} v_i^2 \right )^{1/2} > 0.\]
Further, let $\{(\t{W}_i,v_i)\}_{i \in I}$  be a $(\lambda_1,\lambda_2,\epsilon)$-perturbation
of $\{(W_i,v_i)\}_{i \in I}$ for some $0 \le \l_2 < 1$. Then $\{(\t{W}_i,v_i)\}_{i \in I}$
is a fusion frame with fusion frame bounds
\[
\left [ \frac{ (1- \l_1)\sqrt{C} -
\epsilon \left ( \sum_{i\in I} v_i^2 \right )^{1/2} }{1+\l_2}
\right ]^2
\quad \mbox{and} \quad
\left [ \frac{\sqrt{D}(1+\l_1) +
\epsilon \left ( \sum_{i\in I}v_i^2 \right )^{1/2}}{1-\l_2}\right ]^2.
\]
\end{proposition}

\begin{proof}
We first prove the upper bound. For each $f \in \cH$, we get
\begin{eqnarray*}
\lefteqn{\left[\sum_{i \in I} v_i^2 \|\pi_{\t{W}_i}(f)\|^2 \right]^{1/2}}\\
&\le& \left[ \sum_{i \in I} v_i^2 \left( \norm{\pi_{W_i}(f)} +
\|\pi_{W_i}(f)-\pi_{\t{W}_i}(f)\| \right)^2 \right ]^{1/2}\\
& \le & \left [ \sum_{i \in I} v_i^2 \left( \norm{\pi_{W_i}(f)} +
\l_1 \norm{\pi_{W_i}(f)} + \l_2\|\pi_{\t{W}_i}(f)\| +
\epsilon \|f\| \right)^2 \right]^{1/2}\\
& = & \left [ \sum_{i \in I}  \left( (1+\l_1)v_i \|\pi_{W_i}(f)\| +
\l_2 v_i \|\pi_{\t{W}_i}(f)\| + \epsilon v_i \|f\|
\right)^2 \right ]^{1/2} \\
& \le & \left[ \sum_{i \in I}   (1+\l_1)^2 v_i^2
\|\pi_{W_i}(f)\|^2 \right]^{1/2} +
\left[ \l_2^2 \sum_{i\in I} v_i^2 \|\pi_{\t{W}_i}(f)\|^2 \right]^{1/2}
+ \left( \sum_{i\in I} \epsilon^2 v_i^2 \|f\|^2 \right)^{1/2}\\
&=& (1+ \l_1 ) \left[ \sum_{i\in I}v_i^2 \|\pi_{W_i}(f)\|^2
\right]^{1/2}
+ \l_2 \left[ \sum_{i\in I}v_i^2 \|\pi_{\t{W}_i}(f)\|^2 \right]^{1/2}
+ \epsilon \left( \sum_{i\in I}v_i^2 \right)^{1/2} \hspace*{-0.3cm}\|f\| .
\end{eqnarray*}
Hence,
\[
(1-\l_2) \left [  \sum_{i \in I} v_i^2 \|\pi_{\t{W}_i}(f)\|^2 \right]^{1/2}
\le  (1+\l_1) \left [
\sum_{i \in I} v_i^2\norm{\pi_{W_i}(f)}^2 \right ]^{1/2} +
\epsilon \left ( \sum_{i\in I}v_i^2 \right )^{1/2} \hspace*{-0.3cm}\|f\|,\]
which yields
\[
\left[  \sum_{i \in I} v_i^2 \|\pi_{\t{W}_i}(f)\|^2 \right]^{1/2}
\le
\frac{\sqrt{D}(1+\l_1) +
\epsilon \left( \sum_{i\in I}v_i^2 \right)^{1/2}}{1-\l_2}
 \|f\|.
\]

To prove the lower bound, for all $f \in \cH$ we have
\begin{eqnarray*}
\lefteqn{\left[ \sum_{i \in I} v_i^2 \|\pi_{\t{W}_i}(f)\|^2
\right]^{1/2}}\\
& \ge & \left[ \sum_{i \in I} v_i^2 \left( \norm{\pi_{W_i}(f)} -
\norm{\pi_{W_i}(f)-\pi_{\t{W}_i}(f)} \right)^2 \right ]^{1/2}\\
& \ge & \left[ \sum_{i \in I} v_i^2 \left( \|\pi_{W_i}(f)\| -
\l_1 \norm{\pi_{W_i}(f)} - \l_2\|\pi_{\t{W}_i}(f)\|
- \epsilon \|f\| \right)^2 \right]^{1/2}\\
& = & \left[ \sum_{i \in I}  \left( (1-\l_1)v_i \norm{\pi_{W_i}(f)} -
\l_2 v_i \|\pi_{\t{W}_i}(f)\|
- \epsilon v_i \|f\| \right)^2 \right]^{1/2}\\
& \ge & \left[ \sum_{i \in I} v_i^2  (1-\l_1)^2 \norm{\pi_{W_i}(f)}^2
\right]^{1/2} - \left[ \sum_{i\in I}
 \l_2^2  v_i^2\|\pi_{\t{W}_i}(f)\|^2 \right]^{1/2}
- \left[ \sum_{i\in I} {\epsilon}^2 v_i^2 \|f\|^2 \right]^{1/2}.
\end{eqnarray*}
This implies
\begin{eqnarray*}
{(1+\l_2)\left [ \sum_{i \in I} v_i^2 \|\pi_{\t{W}_i}(f)\|^2
\right ]^{1/2}}
& \ge &   (1-\l_1) \left [ \sum_{i \in J} v_i^2  \norm{\pi_{W_i}(f)}^2
\right ]^{1/2} - \epsilon \left ( \sum_{i\in I} v_i^2 \right )^{1/2}\hspace*{-0.3cm}\|f\| \\
&\ge& (1- \l_1)\sqrt{C}\|f\|
-\epsilon \left ( \sum_{i\in I} v_i^2 \right )^{1/2}\hspace*{-0.3cm}\|f\|\\
&=& \left [ (1- \l_1)\sqrt{C} -
\epsilon \left ( \sum_{i\in I} v_i^2 \right )^{1/2} \right ] \|f\|,
\end{eqnarray*}
which leads to
\[ \left [ \sum_{i \in I} v_i^2 \|\pi_{\t{W}_i}(f)\|^2
\right ]^{1/2} \ge
\frac{\left [ (1- \l_1)\sqrt{C} -
\epsilon \left ( \sum_{i\in I} v_i^2 \right )^{1/2} \right ]}{1+\l_2}
\|f\|.\]
This completes the proof.
\end{proof}

We remark that a different perturbation result for fusion frames can be
derived from \cite[Thm. 3.1]{Sun06} by employing a different definition of
perturbation. However, this would not lead to a result about robustness
of fusion frame systems under disturbances of the local frames, which is what we
are in particular aiming for.


\subsection{Perturbation of the Local Frames}
\label{subsec:pert_frames}

The second fundamental problem with perturbing
fusion frames locally is that a local perturbation cannot ``see'' the
global structure of the fusion frame and therefore cannot adjust for
it.

For the notion of perturbations of sequences we employ the canonical Paley-Wiener-type
definition (compare \cite{CC97}):

\begin{definition}
Let $\{f_i\}_{i \in I}$ and $\{\tilde{f}_i\}_{i \in I}$ be sequences in $\cH$, and
let $0 \le \lambda_1,\lambda_2 < 1$. If
\[
 \norm{\sum_{i \in I} a_i (f_i-\tilde{f}_i)} \le \l_1\norm{\sum_{i \in I} a_i f_i}
 + \l_2 \norm{\sum_{i \in I} a_i \tilde{f}_i}\quad \mbox{for all } \{a_i\}_{i \in I}
 \in \ell^2(I),
\]
then we say that $\{\tilde{f}_i\}_{i \in I}$ is a {\em $(\lambda_1,\lambda_2)$-perturbation}
of $\{f_i\}_{i \in I}$.
\end{definition}

First we derive properties of the relation between the two subspaces spanned by
a sequence and its perturbed version.

\begin{proposition}
\label{prop:f_FF_pert}
Let $\{f_i\}_{i \in I}$ be a frame sequence in $\cH$, and let $0 \le \l_1,\l_2< 1$.
Suppose that $\{\tilde{f}_i\}_{i \in I}$ is a $(\l_1,\l_2)$-perturbation of
$\{f_i\}_{i \in I}$.  Then $\{f_i\}_{i \in I}$ is equivalent to $\{\tilde{f}_i\}_{i \in I}$. In
particular, we have
\[
\frac{1-\l_1}{1+\l_2}\norm{\sum_{i \in I} a_i \tilde{f}_i}
\le \norm{\sum_{i \in I} a_i f_i}
\le \frac{1+\l_2}{1-\l_1}\norm{\sum_{i \in I} a_i \tilde{f}_i} \quad \mbox{for all } \{a_i\}_{i \in I}
 \in \ell^2(I)
\]
and $\dim(\span_{i \in I}\{f_i\}) = \dim(\span_{i \in I}\{\tilde{f}_i\})$.

Set $W=\span_{i \in I}\{f_i\}$ and $\t{W} = \span_{i \in I}\{\tilde{f}_i\}$.
Then
\[
\norm{\pi_W(\pi_{\t{W}}(f))}
\ge \left(\frac{1-\l_1}{1+\l_2}-\l_1\frac{1+\l_2}{1-\l_1} -
\l_2\right) \norm{\pi_{\t{W}}(f)}
\quad \mbox{for all } f \in \cH,
\]
i.e., if $\l_1,\l_2 \le \frac15$, then $\pi_W$ is an isomorphism on Rng $\pi_{\t{W}}$.
\end{proposition}

\begin{proof}
The first part follows from \cite{Chr03}.

Let $S_{\t{\cF}}$ be the frame operator of $\t{\cF}=\{\tilde{f}_i\}_{i \in I}$. For $f \in \cH$,
\begin{eqnarray*}
\norm{\sum_{i \in I} \ip{f}{S_{\t{\cF}}^{-1} \tilde{f}_i} (f_i-\tilde{f}_i)}
& \le & \l_1 \norm{\sum_{i \in I} \ip{f}{S_{\t{\cF}}^{-1} \tilde{f}_i} f_i}
+ \l_2 \norm{\sum_{i \in I} \ip{f}{S_{\t{\cF}}^{-1} \tilde{f}_i} \tilde{f}_i}\\
& \le  & \l_1 \frac{1+\l_2}{1-\l_1} \norm{\sum_{i \in I} \ip{f}{S_{\t{\cF}}^{-1}
\tilde{f}_i} \tilde{f}_i}
+ \l_2 \norm{\sum_{i \in I} \ip{f}{S_{\t{\cF}}^{-1} \tilde{f}_i} \tilde{f}_i}\\
& = & \left(\l_1\frac{1+\l_2}{1-\l_1} + \l_2\right) \norm{\pi_{\t{W}}(f)}.
\end{eqnarray*}
Employing this relation, we obtain
\begin{eqnarray*}
\norm{\pi_W(\pi_{\t{W}}(f))}
& = & \norm{\pi_W\left(\sum_{i \in I} \ip{f}{S_{\t{\cF}}^{-1} \tilde{f}_i}\tilde{f}_i\right)}
\\
& = & \norm{\sum_{i \in I} \ip{f}{S_{\t{\cF}}^{-1} \tilde{f}_i} \pi_W(\tilde{f}_i)}\\
& \ge & \norm{\sum_{i \in I} \ip{f}{S_{\t{\cF}}^{-1} \tilde{f}_i} \pi_W(f_i)}
- \norm{\sum_{i \in I} \ip{f}{S_{\t{\cF}}^{-1} \tilde{f}_i} \pi_W(f_i-\tilde{f}_i)}\\
& = & \norm{\sum_{i \in I} \ip{f}{S_{\t{\cF}}^{-1} \tilde{f}_i} f_i}
- \norm{\sum_{i \in I} \ip{f}{S_{\t{\cF}}^{-1} \tilde{f}_i} \pi_W(f_i-\tilde{f}_i)}\\
& \ge & \frac{1-\l_1}{1+\l_2} \norm{\sum_{i \in I} \ip{f}{S_{\t{\cF}}^{-1} \tilde{f}_i}
\tilde{f}_i}
- \norm{\pi_W\left(\sum_{i \in I} \ip{f}{S_{\t{\cF}}^{-1} \tilde{f}_i}
(f_i-\tilde{f}_i)\right)}\\
& \ge & \frac{1-\l_1}{1+\l_2} \norm{\sum_{i \in I} \ip{f}{S_{\t{\cF}}^{-1} \tilde{f}_i}
\tilde{f}_i}
- \norm{\sum_{i \in I} \ip{f}{S_{\t{\cF}}^{-1} \tilde{f}_i} (f_i-\tilde{f}_i)}\\
& \ge &  \frac{1-\l_1}{1+\l_2} \norm{\pi_{\t{W}}(f)}
- \left(\l_1\frac{1+\l_2}{1-\l_1} + \l_2\right) \norm{\pi_{\t{W}}(f)}\\
& = & \left(\frac{1-\l_1}{1+\l_2}-\l_1\frac{1+\l_2}{1-\l_1} -
\l_2\right)\norm{\pi_{\t{W}}(f)}.
\end{eqnarray*}
It remains to observe that $\l_1,\l_2 \le \frac15$ implies
$\frac{1-\l_1}{1+\l_2}-\l_1\frac{1+\l_2}{1-\l_1} - \l_2> 0$.
\end{proof}

\begin{remark}\label{remaa}
{\rm We wish to mention that in Proposition \ref{prop:f_FF_pert} we did not make
use of the frame bounds of $\{f_i\}_{i\in I}$, but only of the constants
$\l_1,\l_2$ associated with the perturbation. Therefore it follows that our
argument is symmetric in ${\pi}_W$ and ${\pi}_{\t{W}}$ and each permutation
yields the same bounds.}
\end{remark}

The following theorem gives a precise statement of how a perturbation of the local
frames of a fusion frame system -- which certainly results in a perturbation
of the associated fusion frame -- affects its fusion frame bounds,
thereby in particular providing us with conditions under which the subspaces
associated with perturbed local frames still constitute a fusion frame. Thus
fusion frames systems are indeed robust not only against perturbations of
the associated fusion frame (Proposition \ref{propp}), but even against perturbations of the
local frames.

\begin{theorem}
Let $\{(W_i,v_i,\{f_{ij}\}_{j \in J_i})\}_{i\in I}$
be a fusion frame system for $\cH$ with fusion frame bounds $C$, $D$. Choose
$0 \le \l_1,\l_2 < 1$ and $\epsilon >0$ such that
$$
1-\frac{\epsilon^2}{2} = \left(\frac{1-\l_1}{1+\l_2}-
\l_1\frac{1+\l_2}{1-\l_1} -  \l_2\right)
\quad \mbox{and} \quad
\sqrt{C} - \epsilon \left( \sum_{i\in I}v_i^2 \right)^{1/2} >0.
$$
For every $i\in I$, let $\{\tilde{f}_{ij}\}_{j\in J_i}$ be a
$(\l_1,\l_2)$-perturbation of $\{f_{ij}\}_{j\in J_i}$ and let
$\t{W}_i = \span \{\tilde{f}_{ij}\}_{j\in J_i}$.  Then $\{(\t{W}_i,v_i)\}_{i\in I}$
is a fusion frame for $\cH$ with fusion frame bounds
$$
\left[ \sqrt{C} - \epsilon \left( \sum_{i\in I}v_i^2 \right)^{1/2}
\right]^2
\quad \mbox{and} \quad
\left [ \sqrt{D} + \epsilon
\left( \sum_{i\in I}v_i^2 \right)^{1/2} \right]^2.
$$
\end{theorem}

\begin{proof}
Fix $i\in I$. Recalling Proposition \ref{prop:f_FF_pert}, for all $f \in \cH$ we have
\begin{eqnarray*}
 \norm{\pi_{W_i}(f)}^2
   &=&   \norm{\pi_{\t{W}_i}(\pi_{W_i})(f)}^2 + \norm{(I-\pi_{\t{W}_i})\pi_{W_i}(f)}^2\\
   &\ge& \left(1-\frac{\epsilon^2}{2}\right)\norm{\pi_{W_i}(f)}^2 + \norm{(I-\pi_{\t{W}_i})\pi_{W_i}(f)}^2.
\end{eqnarray*}
Hence,
$$
 \norm{(I-\pi_{\t{W}_i})\pi_{W_i}(f)}^2 \le \frac{\epsilon^2}{2}
 \norm{\pi_{W_i}(f)}^2.
$$
Employing Remark \ref{remaa}, Proposition \ref{prop:f_FF_pert} also yields
$$
\norm{(I-\pi_{{W}_i})\pi_{\t{W_i}}(f)}^2 \le
\frac{\epsilon^2}{2} \norm{\pi_{\t{W_i}}(f)}^2.
$$
Collecting the estimates derived above, for any $f\in \cH$ we obtain
\begin{eqnarray*}
\|(\pi_{W_i}-\pi_{\t{W}_i})(f)\|^2 &=& \langle (\pi_{W_i}-\pi_{\t{W}_i})(f),
(\pi_{W_i}-\pi_{\t{W}_i})(f)\rangle \\
&=& \langle (\pi_{W_i}-\pi_{\t{W}_i})^2(f),f\rangle \\
&=& \langle (\pi_{W_i} - \pi_{\t{W}_i}\pi_{W_i} + \pi_{\t{W}_i}
- \pi_{W_i}\pi_{\t{W}_i})(f),f\rangle \\
&\le& \| (I- \pi_{\t{W}_i})(\pi_{W_i}(f)) + (I- \pi_{W_i})(\pi_{\t{W}_i}(f))\|
\|f\|\\
&\le& \|(I- \pi_{\t{W}_i})(\pi_{W_i}(f))\|\|f\| +
\|(I- \pi_{W_i})(\pi_{\t{W}_i}(f))\|\|f\|\\
&\le& \frac{\epsilon^2}{2} \|\pi_{W_i}(f)\|\|f\| +
\frac{\epsilon^2}{2} \|\pi_{\t{W}_i}(f)\| \|f\|\\
&\le& \epsilon^2\|f\|^2.
\end{eqnarray*}
The theorem now follows from Proposition \ref{propp}.
\end{proof}


\section{Discussion}

Fusion frames and fusion frame systems are natural extensions of the theory of
frames.  We have seen pressing needs of such notions and mathematical means in
a variety of applications raging from sensor networks in geophysics, remote sensing
and physiological ear and hearing systems to necessary parallel processing of large frame
systems.   Fusion frames provide a tool for weighted information combination from a set of
overlapping subspaces in a distributed manner.  The fusion process can be either based on
local subspace processing/reconstructions through a fusion frame operator, or based on a global
reconstruction with a distributed dual frame evaluation.  More notably, for fusion frame systems,
the fusion frame operator itself is numerically simple and efficient.  It is a natural
combination of local frame operators, which can be expressed as sums of local frame matrix
representations.  This makes fusion frames more than a mere notion but a practically and
numerically handy tool for distributed processing. Stability is also a substantial feature
of fusion frames, which makes it a robust distributed fusion processing system.

For distributed systems such as general sensor networks, fusion frames are a
tool easy and ready for system modeling and information combination.
For large frame systems requiring parallel processing, fusion frames provide a means to
subdivide the large system into a set of rather flexible and overlapping small
subsystems.   Each subsystem can be processed independently and then combined
coherently.   The computational efficiency is comprehensible.

We envision that applications of fusion frames can reach afar with ample impact.


\section*{Acknowledgments}
Parts of the research for this paper was performed while the second author
was visiting the Department of Mathematics at the University of
Missouri in Columbia. This author thanks this department for its
hospitality and support during this visit. First of all the authors are very
thankful to J. C. Tremain for very helpful suggestions concerning the notion
of perturbation for fusion frames which improved the paper. We are also indebted to
H. B\"olcskei, H. Feichtinger, J. Tanner, and K. Vasudevan for interesting discussions
concerning further applications of fusion frames.

\end{document}